\documentclass[11pt]{article}
\usepackage{latexsym,amssymb,amsmath,enumerate,geometry,float,cite,tikz,setspace,verbatim}
\newtheorem{definition}{Definition}
\newtheorem{theorem}[definition]{Theorem}

\newtheorem{lemma}[definition]{Lemma}

\newtheorem{observation}[definition]{Observation}

\usepackage{lineno}

\begin{document}


\onehalfspace

\title{Strong Equality of Roman and Weak Roman Domination in Trees}

\author{Jos\'{e} D. Alvarado$^1$, Simone Dantas$^1$, Dieter Rautenbach$^2$}

\date{}

\maketitle

\begin{center}
{\small 
$^1$ Instituto de Matem\'{a}tica e Estat\'{i}stica, Universidade Federal Fluminense, Niter\'{o}i, Brazil,
\texttt{josealvarado.mat17@gmail.com, sdantas@im.uff.br}\\[3mm]
$^2$ Institute of Optimization and Operations Research, Ulm University, Ulm, Germany,
\texttt{dieter.rautenbach@uni-ulm.de}
}
\end{center}

\begin{abstract}
We provide a constructive characterization of the trees for which 
the Roman domination number strongly equals the weak Roman domination number,
that is, for which every weak Roman dominating function of minimum weight
is a Roman dominating function.
Our characterization is based on five simple extension operations,
and reveals several structural properties of these trees.
\end{abstract}

{\small 

\medskip


\medskip


}

\section{Introduction}\label{section1}

We consider finite, simple, and undirected graphs, and use standard terminology and notation. 

Let $G$ be a graph, and let $X$ be a subset of the vertex set $V(G)$ of $G$.
For a function $f:V(G)\to \mathbb{R}$, let $f(X)=\sum_{u\in X}f(u)$,
and let the {\it weight} of $f$ be $f(V(G))$.
Furthermore, if $u$ and $v$ are distinct vertices of $G$, then let
$$f_{v\to u}:V(G)\to \mathbb{R}:
x\mapsto
\left\{
\begin{array}{ll}
f(u)+1 &, x=u,\\
f(v)-1 &, x=v,\mbox{ and}\\
f(x) &, x\in V(G)\setminus \{ u,v\}.
\end{array}
\right.
$$
A set $D$ of vertices of $G$ is {\it $X$-dominating} if every vertex in $X\setminus D$ has a neighbor in $D$.
For a positive integer $k$, let $[k]=\{ i\in \mathbb{N}: i\leq k\}$.

Roman domination and weak Roman domination were introduced in \cite{s} and \cite{hehe}, respectively.
For our current purposes, we introduce slightly more general notions.
A {\it Roman dominating function for $(G,X)$}, a $(G,X)${\it-RDF} for short, is a function $f:V(G)\to \{ 0,1,2\}$ 
such that every vertex $u$ in $X$ with $f(u)=0$
has a neighbor $v$ with $f(v)=2$.
The {\it Roman domination number $\gamma_R(G,X)$} of $(G,X)$
is the minimum weight of a $(G,X)$-RDF,
and a $(G,X)$-RDF of weight $\gamma_R(G,X)$ is {\it minimum}.
The {\it Roman domination number $\gamma_R(G)$} of $G$ is $\gamma_R(G,V(G))$.
A {\it weak Roman dominating function for $(G,X)$}, a $(G,X)${\it-WRDF} for short, 
is a function $g:V(G)\to \{ 0,1,2\}$ 
such that every vertex $u$ in $X$ with $g(u)=0$ 
has a neighbor $v$ with $g(v)\geq 1$
such that the set $\{ x\in V(G):g_{v\to u}(x)\geq 1\}$ is $X$-dominating.
The {\it weak Roman domination number $\gamma_r(G,X)$} of $(G,X)$
is the minimum weight of a $(G,X)$-WRDF,
and a $(G,X)$-WRDF of weight $\gamma_r(G,X)$ is {\it minimum}.
The {\it weak Roman domination number $\gamma_r(G)$} of $G$ is $\gamma_r(G,V(G))$.

Since every $(G,X)$-RDF is also a $(G,X)$-WRDF,
we have $\gamma_r(G,X)\leq \gamma_R(G,X)$, 
and, in particular, 
\begin{eqnarray}
\gamma_r(G)\leq \gamma_R(G).\label{e0}
\end{eqnarray}
The motivation for the current work was a problem posed by Chellali et al.~\cite{chh}
who asked for a characterization of the trees that satisfy (\ref{e0}) with equality (cf. Problem 15 in \cite{chh}).
In view of the following result, the extremal graphs for (\ref{e0}) 
do most likely not have a good characterization in general,
which justifies the restriction to trees.

\begin{figure}[H]
\begin{center}
\unitlength 1.3mm 
\linethickness{0.4pt}
\ifx\plotpoint\undefined\newsavebox{\plotpoint}\fi 
\begin{picture}(52,26)(0,0)
\put(5,15){\circle*{1.5}}
\put(24,15){\circle*{1.5}}
\put(43,15){\circle*{1.5}}
\put(10,15){\circle*{1.5}}
\put(29,15){\circle*{1.5}}
\put(48,15){\circle*{1.5}}
\put(5,25){\circle*{1.5}}
\put(24,25){\circle*{1.5}}
\put(43,25){\circle*{1.5}}
\put(10,25){\circle*{1.5}}
\put(29,25){\circle*{1.5}}
\put(48,25){\circle*{1.5}}
\put(5,25){\line(0,-1){10}}
\put(24,25){\line(0,-1){10}}
\put(43,25){\line(0,-1){10}}
\put(5,15){\line(1,2){5}}
\put(24,15){\line(1,2){5}}
\put(43,15){\line(1,2){5}}
\put(10,25){\line(0,-1){10}}
\put(29,25){\line(0,-1){10}}
\put(48,25){\line(0,-1){10}}
\put(10,15){\line(-1,2){5}}
\put(29,15){\line(-1,2){5}}
\put(48,15){\line(-1,2){5}}
\put(5,15){\line(1,0){5}}
\put(24,15){\line(1,0){5}}
\put(43,15){\line(1,0){5}}
\put(1,15){\makebox(0,0)[cc]{$x_1$}}
\put(20,15){\makebox(0,0)[cc]{$x_2$}}
\put(39,15){\makebox(0,0)[cc]{$x_3$}}
\put(14,15){\makebox(0,0)[cc]{$\bar{x}_1$}}
\put(33,15){\makebox(0,0)[cc]{$\bar{x}_2$}}
\put(52,15){\makebox(0,0)[cc]{$\bar{x}_3$}}
\put(17,5){\circle*{1.5}}
\put(37,5){\circle*{1.5}}
\put(17,1){\makebox(0,0)[cc]{$c_1$}}
\put(5,15){\line(6,-5){12}}
\multiput(17,5)(.033653846,.048076923){208}{\line(0,1){.048076923}}
\multiput(48,15)(-.1043771044,-.0336700337){297}{\line(-1,0){.1043771044}}
\put(37,1){\makebox(0,0)[cc]{$c_2$}}
\multiput(10,15)(.0909090909,-.0336700337){297}{\line(1,0){.0909090909}}
\multiput(37,5)(-.0437710438,.0336700337){297}{\line(-1,0){.0437710438}}
\multiput(48,15)(-.037037037,-.0336700337){297}{\line(-1,0){.037037037}}
\end{picture}
\end{center}
\caption{The graph $G$ for the two clauses 
$C_1=x_1\vee x_2\vee \bar{x}_3$ and $C_2=\bar{x}_1\vee x_2\vee \bar{x}_3$ over the three boolean variables $x_1$, $x_2$, and $x_3$. 
}\label{fig0}
\end{figure}
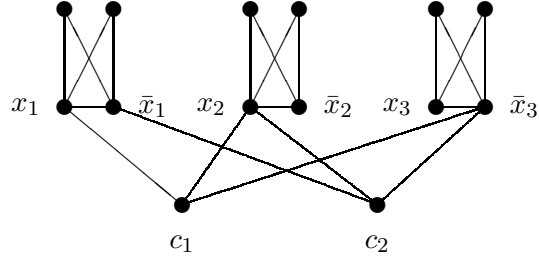

\begin{theorem}\label{theorem1}
For a given graph $G$, it is NP-hard to decide whether $\gamma_r(G)=\gamma_R(G)$.
\end{theorem}
{\it Proof:} We describe a polynomial reduction from {\sc 3Sat}. Therefore, let $F$ be a {\sc 3Sat} instance with 
clauses $C_1,\ldots,C_m$ over the boolean variables $x_1,\ldots,x_n$. We construct a graph $G$ whose order is polynomially bounded in terms of $n$ and $m$ such that 
$F$ is satisfiable if and only if $\gamma_r(G)=\gamma_R(G)$.
Therefore, for every boolean variable $x_i$, 
create a copy $G(x_i)$ of $K_4-e$ and denote the two vertices of degree $3$ in $G(x_i)$ by $x_i$ and $\bar{x}_i$.
For every clause $C_j$, create a vertex $c_j$.
For every literal $x\in \{ x_i,\bar{x}_i\}$ and every clause $C_j$ such that $x$ appears in $C_j$, connect the vertex denoted $x$ in $G(x_i)$ with $c_j$ by an edge.
See Figure \ref{fig0} for an example of the construction.

\noindent Clearly, for every $(G,V(G))$-WRDF $g$ and every $i\in [n]$, 
we have $g(V(G(x_i)))\geq 2$, which implies 
$2n\leq \gamma_r(G)\leq \gamma_R(G)$.
Since 
$$g:V(G)\to \{ 0,1,2\}:x\mapsto
\left\{
\begin{array}{ll}
1 &, x\in \{ x_i:i\in [n]\}\cup \{ \bar{x}_i:i\in [n]\}\mbox{, and }\\
0 &, x\in V(G)\setminus \left(\{ x_i:i\in [n]\}\cup \{ \bar{x}_i:i\in [n]\}\right)
\end{array}
\right.
$$
is a $(G,V(G))$-WRDF, we have $\gamma_r(G)=2n$.
Furthermore, $\gamma_R(G)=2n$ holds if and only if there is a $(G,V(G))$-RDF $f$ such that 
for every $i\in [n]$, 
$f$ assigns the value $2$ to either $x_i$ or $\bar{x}_i$,
and to every other vertex, $f$ assigns the value $0$.
Since such a $(G,V(G))$-RDF indicates a satisfying truth assignment for $F$, 
and, conversely, 
a satisfying truth assignment for $F$ leads to such a $(G,V(G))$-RDF, 
we obtain that 
$\gamma_r(G)=\gamma_R(G)$
if and only if $F$ is satisfiable. $\Box$

\bigskip

\noindent A typical solution for the problem posed by Chellali et al.~\cite{chh}
would be a so-called {\it constructive characterization}, 
that is, a recursive constructive description of the set of all extremal trees for (\ref{e0}). 
There are many examples of such characterizations in the literature \cite{cfmp,dghm,hh,h}.
Usually, they involve some few small extremal trees together 
with a small set of simple extension operations 
that are applied recursively 
in order to create all larger extremal trees. 
Sometimes additional information, such as certain labels or suitable subsets, 
has to be maintained in order to apply the extension operations properly. 
The Roman domination number as well as the weak Roman domination number of a given tree
can be calculated by simple linear time algorithms based on standard approaches \cite{hehe}.
This implies that the extremal trees for (\ref{e0}) can easily be recognized in linear time,
and a constructive characterization of these trees 
would only be beneficial 
if it reveals interesting structural properties and/or is considerably simpler than the two linear time algorithms.
We did not arrive at a completely satisfactory solution of the problem posed by Chellali et al.~\cite{chh},
because all our constructive characterizations 
were essentially equivalent to implicit executions of the two linear time algorithms.
Therefore, we turn to a variation of the posed problem
based on the concept of {\it strong equality}, 
which was first introduced by Haynes and Slater in \cite{hs}.

The Roman domination number of a graph $G$ {\it strongly equals} the weak Roman domination number of $G$
if every minimum $(G,V(G))$-WRDF is a $(G,V(G))$-RDF.
Since the Roman domination number of $G$ equals the weak Roman domination number of $G$
if some - and not necessarily all - minimum $(G,V(G))$-WRDF is a $(G,V(G))$-RDF,
strong equality implies equality.
Our main result presented in the next section is a constructive characterization, 
based on five simple extension operations, 
of the trees for which the Roman domination number strongly equals the weak Roman domination number.
Further examples of characterizations of strong equalities can be found in \cite{cr,hhs,rl}.
In a concluding section we discuss a possible constructive characterization of the extremal trees for (\ref{e0})
and its weaknesses.

\section{Constructive characterization of strong equality}\label{section2}

Instead of just trees our construction involves slightly more general objects,
which are trees together with two suitable vertex subsets.
Therefore, let ${\cal S}$ be the set of all triples $(T,X,Y)$ with the following properties:
\begin{itemize}
\item $T$ is a tree, and $X$ and $Y$ are sets of vertices of $T$.
\item Every minimum $(T,X)$-WRDF is a $(T,X)$-RDF.
\item $Y$ is the set of all vertices $u$ of $T$ 
for which there is some minimum $(T,X)$-WRDF $g$ 
such that 
\begin{itemize}
\item either $g(u)\geq 1$, 
\item or $g(u)=0$ and $u$ has a neighbor $v$ with $g(v)\geq 1$
such that the set $\{ x\in V(T):g_{v\to u}(x)\geq 1\}$ is $X$-dominating.
\end{itemize}
\end{itemize}
The definition immediately implies the following observation.

\begin{observation}
The Roman domination number strongly equals the weak Roman domination number for some tree $T$
if and only if ${\cal S}$ contains the triple $(T,V(T),V(T))$.
\end{observation}
The following lemma collects some properties of the elements of ${\cal S}$.

\begin{lemma}\label{lemma1}
For $(T,X,Y)\in {\cal S}$, the following statements hold.
\begin{enumerate}[(i)] 
\item $T$ and $X$ uniquely determine $Y$, and $X\subseteq Y$.
\item Either $|X|=0$, or $|X|=|Y|=|V(T)|=1$, or $|X|\geq 3$.
\item If $g$ is a minimum $(T,X)$-WRDF and $g(u^*)=1$ for some vertex $u^*$ of $T$,
then $V(T)=X=\{ u^*\}$.
\item If  
either $T$ has order $1$ and $X=\emptyset$, 
or $T$ has order at least $2$, then some vertex $u$ of $T$ belongs to $Y$
if and only if there is some minimum $(T,X)$-WRDF $g$ such that $N_T[u]$ contains a vertex $v$ with $g(v)=2$.
\end{enumerate}
\end{lemma}
{\it Proof:} (i) This follows immediately from the definition of ${\cal S}$.

\bigskip

\noindent (ii) Let $T$ be a tree, and let $X$ be a set of two vertices of $T$, say $x_1$ and $x_2$.
If $x_1$ and $x_2$ are adjacent,  
then 
$$g:V(T)\to \{ 0,1,2\}:
x\mapsto
\left\{
\begin{array}{ll}
1 &, x=x_1,\mbox{ and}\\
0 &, x\in V(T)\setminus \{ x_1\},
\end{array}
\right.
$$
and, if $x_1$ and $x_2$ are not adjacent but $y$ is a neighbor of $x_2$,  
then 
$$g:V(T)\to \{ 0,1,2\}:
x\mapsto
\left\{
\begin{array}{ll}
1 &, x\in \{ x_1,y\},\mbox{ and}\\
0 &, x\in V(T)\setminus \{ x_1,y\}
\end{array}
\right.
$$
is a minimum $(T,X)$-WRDF that is not a $(T,X)$-RDF.
This implies that ${\cal S}$ contains no triple $(T,X,Y)$ with $|X|=2$.
Similarly, it follows that ${\cal S}$ contains no triple $(T,X,Y)$ with $|X|=1$ 
where $T$ has order at least $2$.

\bigskip

\noindent (iii) If $T$ has exactly one vertex, then the statement is trivial.
Hence, we may assume for a contradiction, that $T$ has order at least $2$,
$g$ is a minimum $(T,X)$-WRDF, and $g(u^*)=1$ for some $u^*\in V(T)$.
Since every vertex $u$ of $T$ with $g(u)=0$ has a neighbor $v$ with $g(v)=2$,
and the function
$$x\mapsto
\left\{
\begin{array}{ll}
0 &, x=u^*,\mbox{ and}\\
g(x) &, x\in V(T)\setminus \{ u^*\}
\end{array}
\right.
$$
is not a $(T,X)$-WRDF,
we obtain that $u^*\in X$, and that $u^*$ has no neighbor $v$ with $g(v)\geq 1$.
Now, if $v^*$ is any neighbor of $u^*$, then
the function
$$x\mapsto
\left\{
\begin{array}{ll}
0 &, x=u^*,\\
1 &, x=v^*,\mbox{ and}\\
g(x) &, x\in V(T)\setminus \{ u^*,v^*\}
\end{array}
\right.
$$
is a minimum $(T,X)$-WRDF that is not a $(T,X)$-RDF,
which is a contradiction.

\bigskip

\noindent (iv) The ``if'' part of the statement is trivial, and the ``only if'' part follows from (ii). 
$\Box$

\bigskip

\noindent The following two lemmas capture the reduction operations for ${\cal S}$.
While the first lemma is the key result for our constructive characterization,
the second lemma allows to decompose the reduction described in the first lemma
into more elementary reductions, removing only between one and four vertices at a time.

\begin{figure}[H]
\begin{center}
\unitlength 1mm 
\linethickness{0.4pt}
\ifx\plotpoint\undefined\newsavebox{\plotpoint}\fi 
\begin{picture}(70,51)(0,0)
\put(5,15){\circle*{2}}
\put(17,15){\circle*{2}}
\put(41,15){\circle*{2}}
\put(65,15){\circle*{2}}
\put(5,6){\makebox(0,0)[cc]{$W_1$}}
\put(5,12){\makebox(0,0)[cc]{$w_1$}}
\put(17,6){\makebox(0,0)[cc]{$W_2$}}
\put(17,12){\makebox(0,0)[cc]{$w_2$}}
\put(41,6){\makebox(0,0)[cc]{$W_\ell$}}
\put(41,12){\makebox(0,0)[cc]{$w_\ell$}}
\put(65,6){\makebox(0,0)[cc]{$W_k$}}
\put(65,12){\makebox(0,0)[cc]{$w_k$}}
\put(29,11){\makebox(0,0)[cc]{$\cdots$}}
\put(53,11){\makebox(0,0)[cc]{$\cdots$}}
\put(35,25){\circle*{2}}
\put(35,35){\circle*{2}}
\put(65,15){\line(-3,1){30}}
\put(35,25){\line(3,-5){6}}
\put(35,35){\line(0,-1){10}}
\multiput(35,25)(-.0606060606,-.0336700337){297}{\line(-1,0){.0606060606}}
\put(35,25){\line(-3,-1){30}}
\put(38,27){\makebox(0,0)[cc]{$v$}}
\put(35,46){\makebox(0,0)[cc]{$T'$}}
\put(35,38){\makebox(0,0)[cc]{$u$}}
\put(5,10){\oval(10,16)[]}
\put(17,10){\oval(10,16)[]}
\put(41,10){\oval(10,16)[]}
\put(65,10){\oval(10,16)[]}
\put(35,41.5){\oval(12,19)[]}
\end{picture}
\end{center}
\caption{The configuration as in Lemma \ref{lemma2}.}\label{fig2}
\end{figure}
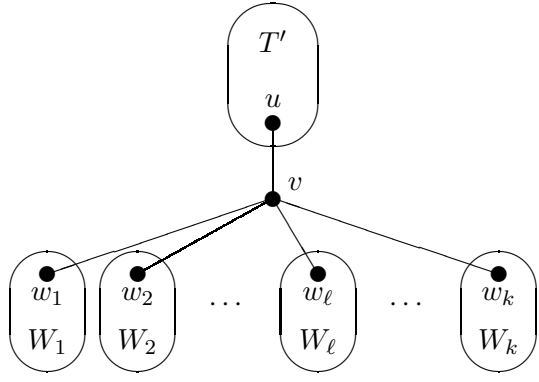

\begin{lemma}\label{lemma2}
Let $T$ be a tree, and let $X$ and $Y$ be sets of vertices of $T$.
Let $v$ be a vertex of $T$, and let $N_T(v)=\{ u,w_1,\ldots,w_k\}$.
For $i\in [k]$, let $W_i$ be the vertex set of the component of $T-v$ that contains $w_i$,
and let $T'$ be the component of $T-v$ that contains $u$ (cf.~Figure \ref{fig2}). 
For some $\ell\in [k]$, let $W_i\cap X=\{ w_i\}$ for $i\in [\ell]$,
and let $W_j\cap X=\emptyset$ for $j\in [k]\setminus [\ell]$.

$(T,X,Y)\in {\cal S}$ if and only if
\begin{enumerate}[(a)]
\item either 
\begin{enumerate}[(i)]
\item $\ell=2$, $u\in X$, $u,v\in Y$, and $W_i\cap Y=\{ w_i\}$ for $i\in [k]$.
\item $(T',X',Y')\in {\cal S}$ for 
$X'=X\setminus \{ u,v,w_1,\ldots,w_\ell\}$ and 
$Y'=(Y\cap V(T'))\setminus \{ u\}$.
\end{enumerate}
\item or 
\begin{enumerate}[(i)]
\item $\ell\geq 3$, $u,v\in Y$, and $W_i\cap Y=\{ w_i\}$ for $i\in [k]$.
\item $(T',X',Y')\in {\cal S}$ for $X'=X\setminus \{ u,v,w_1,\ldots,w_\ell\}$ 
and some $Y'$ with 
$$(Y\cap V(T'))\setminus \{ u\}\subseteq Y'\subseteq Y\cap V(T').$$
\end{enumerate}
\end{enumerate}
\end{lemma}
{\it Proof:} First, let $\ell=1$.
Let $g$ be a minimum $(T,X)$-WRDF.
Clearly, $1\leq g(V(T)\setminus V(T'))\leq 2$.
If $g(V(T)\setminus V(T'))=1$, then 
$$\tilde{g}:V(T)\to \{ 0,1,2\}:
x\mapsto
\left\{
\begin{array}{ll}
1 &, x=v,\\
0 &, x\in V(T)\setminus (V(T')\cup \{ v\}),\mbox{ and}\\
g(x) &, x\in V(T')
\end{array}
\right.
$$
is a minimum $(T,X)$-WRDF.
Since $\tilde{g}$ is not a $(T,X)$-RDF, we have $(T,X,Y)\not\in {\cal S}$ in this case.
If $g(V(T)\setminus V(T'))=2$, then
$$\tilde{g}:V(T)\to \{ 0,1,2\}:
x\mapsto
\left\{
\begin{array}{ll}
1 &, x=v,\\
\min\{ 2,g(u)+1\} &, x=u,\\
0 &, x\in V(T)\setminus (V(T')\cup \{ v\}),\mbox{ and}\\
g(x) &, x\in V(T')\setminus \{ u\}
\end{array}
\right.
$$
is a minimum $(T,X)$-WRDF.
Since $\tilde{g}$ is not a $(T,X)$-RDF, we have $(T,X,Y)\not\in {\cal S}$ also in this case.
Therefore, $\ell\geq 2$ is a necessary condition for $(T,X,Y)\in {\cal S}$.

Let $\ell\geq 2$.

Let $X'=X\setminus \{ u,v,w_1,\ldots,w_\ell\}$.
Note that $X'=X\cap (V(T')\setminus \{ u\})$.

Since $\ell\geq 2$, we obtain, 
for every minimum $(T,X)$-WRDF $g$, that 
$g(V(T)\setminus V(T'))=2$, 
and that the restriction $g\mid_{V(T')}$ of $g$ to $V(T)$ is a minimum $(T',X')$-WRDF,
that is, $\gamma_r(T,X)=\gamma_r(T',X')+2$.

We consider two cases according to the value of $\ell$.

\bigskip

\noindent {\bf Case 1} $\ell=2$.

\bigskip

\noindent First, we prove the necessity, that is, 
we assume that $(T,X,Y)\in {\cal S}$ holds, 
and show that (a)(i) and (a)(ii) hold.
If $u\not\in X$, and $g'$ is a minimum $(T',X')$-WRDF, then
$$g:V(T)\to \{ 0,1,2\}:
x\mapsto
\left\{
\begin{array}{ll}
1 &, x\in \{ v,w_1\},\\
0 &, x\in V(T)\setminus (V(T')\cup \{ v,w_1\}),\mbox{ and}\\
g'(x) &, x\in V(T')
\end{array}
\right.
$$
is a minimum $(T,X)$-WRDF.
Since $g$ is not a $(T,X)$-RDF, we obtain a contradiction.
Hence, $u\in X$ holds.
Since $(T,X,Y)\in {\cal S}$, every minimum $(T,X)$-WRDF $g$ satisfies
$$
g(x)=
\left\{
\begin{array}{ll}
2 &, x=v,\mbox{ and}\\
0 &, x\in V(T)\setminus (V(T')\cup \{ v\}),
\end{array}
\right.
$$
which implies that (a)(i) holds.
If there is no set $Y'$ such that $(T',X',Y')\in {\cal S}$, 
then there is some minimum $(T',X')$-WRDF $g'$ that is not a $(T',X')$-RDF,
that is, there is some $u'\in X'$ with $g'(u')=0$ 
that has no neighbor $v'$ in $T'$ with $g'(v')=2$.
Now
$$g:V(T)\to \{ 0,1,2\}:
x\mapsto
\left\{
\begin{array}{ll}
2 &, x=v,\\
0 &, x\in V(T)\setminus (V(T')\cup \{ v\}),\mbox{ and}\\
g'(x) &, x\in V(T')
\end{array}
\right.
$$
is a minimum $(T,X)$-WRDF.
Since $u'\not=u$, the function $g$ is not a $(T,X)$-RDF, which is a contradiction.
Hence, $(T',X',Y')\in {\cal S}$ for some set $Y'$.
If $u\in Y'$, then, by Lemma \ref{lemma1}(iv), 
there is some minimum $(T',X')$-WRDF $g'$ 
such that $g'(v')=2$ for some $v'\in N_{T'}[u]$.
Now
$$g:V(T)\to \{ 0,1,2\}:
x\mapsto
\left\{
\begin{array}{ll}
1 &, x\in \{ v,w_1\},\\
0 &, x\in V(T)\setminus (V(T')\cup \{ v,w_1\}),\mbox{ and}\\
g'(x) &, x\in V(T')
\end{array}
\right.
$$
is a minimum $(T,X)$-WRDF.
Since $g$ is not a $(T,X)$-RDF, we obtain a contradiction.
Hence, $u\not\in Y'$.
If $u'\in Y'\setminus ((Y\cap V(T'))\setminus \{ u\})$, 
then, by Lemma \ref{lemma1}(iv), 
there is some minimum $(T',X')$-WRDF $g'$
such that $g'(v')=2$ for some $v'\in N_{T'}[u']$.
Now
$$g:V(T)\to \{ 0,1,2\}:
x\mapsto
\left\{
\begin{array}{ll}
2 &, x=v,\\
0 &, x\in V(T)\setminus (V(T')\cup \{ v\}),\mbox{ and}\\
g'(x) &, x\in V(T')
\end{array}
\right.
$$
is a minimum $(T,X)$-WRDF
such that $g(v')=2$ and $v'\in N_T[u']$,
which implies the contradiction $u'\in Y$.
Hence,  $Y'\subseteq ((Y\cap V(T'))\setminus \{ u\})$.
If $u'\in ((Y\cap V(T'))\setminus \{ u\})\setminus Y'$,
then there is some minimum $(T,X)$-WRDF $g$
such that $g(v')=2$ for some $v'\in N_T[u']$.
Now $g\mid_{V(T')}$ is a minimum $(T',X')$-WRDF $g'$
such that $g'(v')=2$ and $v'\in N_{T'}[u']$,
which implies the contradiction $u'\in Y'$.
Altogether, we obtain $Y'=((Y\cap V(T'))\setminus \{ u\})$,
that is, (a)(ii) holds, which completes the proof of the necessity.

\bigskip

\noindent We proceed to the proof of the sufficiency, that is, 
we assume that (a)(i) and (a)(ii) hold, and show that $(T,X,Y)\in {\cal S}$ holds.
By (a)(i) and (a)(ii), we have 
$Y\setminus (V(T')\setminus \{ u\})=\{ u,v,w_1,\ldots,w_k\}$
and $Y=Y'\cup \{ u,v,w_1,\ldots,w_k\}$.
Let $g$ be a minimum $(T,X)$-WRDF. 
Let $g'$ be $g\mid_{V(T')}$.
Recall that $g(V(T)\setminus V(T'))=2$, and that $g'$ is a minimum $(T',X')$-WRDF.
By (a)(ii), the function $g'$ is a $(T',X')$-RDF.
Furthermore, also by (a)(ii), we have $u\not\in Y'$, 
which implies that $g'(u)=0$, and that $u$ has no neighbor $v'$ in $T'$ with $g'(v')\geq 1$
such that $\{ x\in V(T'):g'_{v'\to u}(x)\geq 1\}$ is $X'$-dominating.
Since every vertex $u''$ in $X'$ with $g'(u'')=0$ has a neighbor $v''$ with $g'(v'')=2$,
this implies that $g'(v')=0$ for every $v'\in N_{T'}[u]$.
Since $u,w_1,w_2\in X$, this implies
$$
g(x)=
\left\{
\begin{array}{ll}
2 &, x=v,\mbox{ and}\\
0 &, x\in V(T)\setminus (V(T')\cup \{ v\}),
\end{array}
\right.
$$
which implies that $g$ is a $(T,X)$-RDF.
Hence, $(T,X,\tilde{Y})$ for some set $\tilde{Y}$ 
with $\tilde{Y}\setminus (V(T')\setminus \{ u\})=\{ u,v,w_1,\ldots,w_k\}=Y\setminus (V(T')\setminus \{ u\})$.
It remains to show that $\tilde{Y}=Y$.
If $u'\in \tilde{Y}\setminus Y$,
then $u'\in V(T')\setminus \{ u\}$, and, by Lemma \ref{lemma1}(iv),
there is a minimum $(T,X)$-WRDF $g$ 
such that $g(v')=2$ for some $v'\in N_T[u']$.
Now $g\mid_{V(T')}$ is a minimum $(T',X')$-WRDF $g'$
such that $g'(v')=2$ and $v'\in N_{T'}[u']$,
which implies the contradiction $u'\in Y$.
If $u'\in Y\setminus \tilde{Y}$,
then $u'\in V(T')\setminus \{ u\}$, and, by Lemma \ref{lemma1}(iv),
there is a minimum $(T',X')$-WRDF $g'$ 
such that $g(v')=2$ for some $v'\in N_{T'}[u']$.
Now 
$$g:V(T)\to \{ 0,1,2\}:
x\mapsto
\left\{
\begin{array}{ll}
2 &, x=v,\\
0 &, x\in V(T)\setminus (V(T')\cup \{ v\}),\mbox{ and}\\
g'(x) &, x\in V(T')
\end{array}
\right.
$$
is a minimum $(T,X)$-WRDF $g$
such that $g(v')=2$ and $v'\in N_T[u']$,
which implies the contradiction $u'\in \tilde{Y}$.
Altogether, we obtain $\tilde{Y}=Y$,
which completes the proof in this case.

\bigskip

\noindent {\bf Case 2} $\ell\geq 3$.

\bigskip

\noindent Since the proof in this case is similar to - and simpler than - the proof in Case 1, 
we leave some details to the reader.

First, we prove the necessity, and assume that $(T,X,Y)\in {\cal S}$ holds.
Since $\ell\geq 3$, we obtain for every minimum $(T,X)$-WRDF $g$, that 
$$
g(x)=
\left\{
\begin{array}{ll}
2 &, x=v,\mbox{ and}\\
0 &, x\in V(T)\setminus (V(T')\cup \{ v\}),
\end{array}
\right.
$$
which implies (b)(i).
Exactly as in the proof for Case 1, we obtain that $(T',X',Y')\in {\cal S}$ for some set $Y'$.
If $u'\in Y'$, 
then similar arguments as in Case 1 imply $u'\in Y\cap V(T')$,
which implies $Y'\subseteq Y\cap V(T')$.
If $u'\in (Y\cap V(T'))\setminus \{ u\}$, 
then similar arguments as in Case 1 imply $u'\in Y'$,
which implies $(Y\cap V(T'))\setminus \{ u\}\subseteq Y'$.
Altogether, we obtain that (b)(ii) holds,
which completes the proof of the necessity.

\bigskip

\noindent Next, we prove the sufficiency, and assume that (b)(i) and (b)(ii) hold.
If $g$ is a minimum $(T,X)$-WRDF, and
$g'$ is $g\mid_{V(T')}$,
then 
$$
g(x)=
\left\{
\begin{array}{ll}
2 &, x=v,\mbox{ and}\\
0 &, x\in V(T)\setminus (V(T')\cup \{ v\}),
\end{array}
\right.
$$
and $g'$ is a minimum $(T',X')$-WRDF.
By (b)(ii), the function $g'$ is a $(T',X')$-RDF, 
and hence, $g$ is a $(T,X)$-RDF.
This implies that $(T,X,\tilde{Y})$ for some set $\tilde{Y}$ 
with $\tilde{Y}\setminus (V(T')\setminus \{ u\})=\{ u,v,w_1,\ldots,w_k\}=Y\setminus (V(T')\setminus \{ u\})$.
Note that, by (b)(ii), we have $Y\cap (V(T')\setminus \{ u\})=Y'\setminus \{ u\}$.
Therefore, in order to show $\tilde{Y}=Y$,
it remains to show that $\tilde{Y}\cap (V(T')\setminus \{ u\})=Y'\setminus \{ u\}$,
which can be done using similar arguments as in Case 1.
This completes the proof. $\Box$

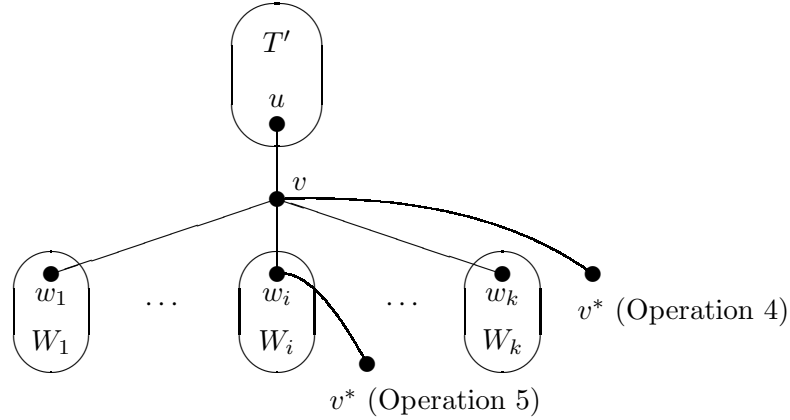
\begin{figure}[H]
\begin{center}
\unitlength 1mm 
\linethickness{0.4pt}
\ifx\plotpoint\undefined\newsavebox{\plotpoint}\fi 
\begin{picture}(78,51)(0,0)
\put(5,15){\circle*{2}}
\put(35,15){\circle*{2}}
\put(65,15){\circle*{2}}
\put(5,6){\makebox(0,0)[cc]{$W_1$}}
\put(5,12){\makebox(0,0)[cc]{$w_1$}}
\put(35,6){\makebox(0,0)[cc]{$W_i$}}
\put(35,12){\makebox(0,0)[cc]{$w_i$}}
\put(65,6){\makebox(0,0)[cc]{$W_k$}}
\put(65,12){\makebox(0,0)[cc]{$w_k$}}
\put(35,25){\circle*{2}}
\put(35,35){\circle*{2}}
\put(65,15){\line(-3,1){30}}
\put(35,35){\line(0,-1){10}}
\put(35,25){\line(-3,-1){30}}
\put(38,27){\makebox(0,0)[cc]{$v$}}
\put(35,46){\makebox(0,0)[cc]{$T'$}}
\put(35,38){\makebox(0,0)[cc]{$u$}}
\put(5,10){\oval(10,16)[]}
\put(35,10){\oval(10,16)[]}
\put(65,10){\oval(10,16)[]}
\put(35,41.5){\oval(12,19)[]}
\put(77,15){\circle*{2}}
\put(75,10){\makebox(0,0)[lc]{$v^*$ (Operation 4)}}
\put(42,-2){\makebox(0,0)[lc]{$v^*$ (Operation 5)}}
\qbezier(35,25)(62,26)(77,15)
\qbezier(35,15)(40,16)(47,3)
\put(47,3){\circle*{2}}
\put(52,11){\makebox(0,0)[cc]{$\cdots$}}
\put(20,11){\makebox(0,0)[cc]{$\cdots$}}
\put(35,25){\line(0,-1){10}}
\end{picture}
\end{center}
\caption{The configuration as in Lemma \ref{lemma3}(2),
indicating the two possible types of neighbors of the new endvertex $v^*$.}\label{fig3}
\end{figure}

\begin{lemma}\label{lemma3}
Let $T$ be a tree, and let $X$ and $Y$ be sets of vertices of $T$.
Let $u^*$ be a vertex of $T$, 
and let $T^*$ arise from $T$ by adding a new vertex $v^*$, and a new edge $u^*v^*$,
that is, $v^*$ is an endvertex of $T^*$, and $T=T^*-v^*$.

If one of the following two conditions (1) and (2) is satisfied, then 
$$(T^*,X^*,Y^*)\in {\cal S}\mbox{ if and only if }(T,X,Y)\in {\cal S}.$$
\begin{enumerate}[(1)]
\item $u^*\not\in Y$, $X^*=X$, and $Y^*=Y$.
\item $T$, $u$, $v$, $w_1$, $\ldots$, $w_k$, $W_1$, $\ldots$, $W_k$, and $\ell$ 
are as in the statement of Lemma \ref{lemma2},\\
either (a)(i) or (b)(i) is satisfied,
and
\begin{enumerate}[(2.1)]
\item either $u^*=v$, $X\subseteq X^*\subseteq X\cup \{ v^*\}$, and $Y^*=Y\cup \{ v^*\}$.
\item or $u^*\in \{ w_1,\ldots,w_k\}$, $X^*=X$, and $Y^*=Y$.
\end{enumerate}
(See Figure \ref{fig3} for an illustration.)
\end{enumerate}
\end{lemma}
{\it Proof:} First, we assume that (1) holds, that is, $u^*\not\in Y$, $X^*=X$, and $Y^*=Y$.
In this case, 
for every minimum $(T^*,X^*)$-WRDF $g^*$,
we have $g^*(v^*)=0$, and the restriction $g^*\mid_{V(T)}$ of $g^*$ to $V(T)$ is a minimum $(T,X)$-WRDF.
Conversely, 
for every minimum $(T,X)$-WRDF $g$,
the function 
$$
x\mapsto
\left\{
\begin{array}{ll}
0 &, x=v^*,\mbox{ and}\\
g(x) &, x\in V(T)
\end{array}
\right.
$$
is a minimum $(T^*,X^*)$-WRDF.
Therefore, 
if $(T,X,Y)\in {\cal S}$, and $g^*$ is some minimum $(T^*,X^*)$-WRDF,
then $g^*\mid_{V(T)}$ is a minimum $(T,X)$-WRDF, which is a $(T,X)$-RDF, and, 
hence, $g^*$ is a $(T^*,X^*)$-RDF,
that is, $(T^*,X^*,Y^*)\in {\cal S}$.
Conversely, 
if $(T^*,X^*,Y^*)\in {\cal S}$, and $g$ is some minimum $(T,X)$-WRDF,
then extending $g$ by $0$ on $v^*$ yields a minimum $(T^*,X^*)$-WRDF,
which is a $(T^*,X^*)$-RDF, and, hence, $g$ is a $(T,X)$-RDF,
that is, $(T,X,Y)\in {\cal S}$.
Altogether, $(T^*,X^*,Y^*)\in {\cal S}$ if and only if $(T,X,Y)\in {\cal S}$.

Now, let $T$, $u$, $v$, $w_1$, $\ldots$, $w_k$, $W_1$, $\ldots$, $W_k$, and $\ell$ 
be as in the statement of Lemma \ref{lemma2}, and let either (a)(i) or (b)(i) be satisfied.
Note that if (2.1) holds, that is,
$u^*=v$, $X\subseteq X^*\subseteq X\cup \{ v^*\}$, and $Y^*=Y\cup \{ v^*\}$,
then attaching the new endvertex $v^*$ to $T$ at $v$
is actually equivalent to increasing $k$ and possibly also $\ell$ by $1$, 
that is, $v^*$ plays the role of some additional neighbor of $v$ next to $w_1,\ldots,w_k$.
In fact, if $X^*=X$, then $\ell$ remains unchanged, and 
if $X^*=X\cup \{ v^*\}$, then $\ell$ is increased by $1$.
Similarly, if (2.2) holds,
that is, $u^*\in \{ w_1,\ldots,w_k\}$, $X^*=X$, and $Y^*=Y$,
then $v^*$ is just added as one further vertex to one of the sets $W_i$.
Therefore, both cases describe an extension of the structure of $T-V(T')$ 
specified in Lemma \ref{lemma2} by exactly one vertex.
See Figure \ref{fig3} for an illustration.
Now, the proof can be completed by some applications of Lemma \ref{lemma2}.
If $(T,X,Y)\in {\cal S}$, 
then, 
by the forward implication of Lemma \ref{lemma2}, 
there is some $(T',X',Y')\in {\cal S}$
satisfying Lemma \ref{lemma2}(a) or Lemma \ref{lemma2}(b),
and, by the backward implication of Lemma \ref{lemma2}, $(T^*,X^*,Y^*)\in {\cal S}$.
Conversely, 
if $(T^*,X^*,Y^*)\in {\cal S}$, 
then, 
by the forward implication of Lemma \ref{lemma2}, 
there is some $(T',X',Y')\in {\cal S}$
satisfying Lemma \ref{lemma2}(a) or Lemma \ref{lemma2}(b) where $(T^*,X^*,Y^*)$ replaces $(T,X,Y)$,
and, by the backward implication of Lemma \ref{lemma2}, $(T,X,Y)\in {\cal S}$.
Altogether, $(T^*,X^*,Y^*)\in {\cal S}$ if and only if $(T,X,Y)\in {\cal S}$,
which completes the proof. 
$\Box$

\bigskip

\noindent We are now in a position to describe the five extension operations.
Therefore, let $(T,X,Y)$ be such that $T$ is a tree,
$X$ and $Y$ are sets of vertices of $T$,
and $X\subseteq Y$.
\begin{itemize}
\item {\bf Operation 1}

$(T^+,X^+,Y^+)$ arises by applying Operation 1 to $(T,X,Y)$ if 
\begin{itemize}
\item there is some vertex $u$ of $T$ with $u\not\in Y$,
\item $T^+$ arises from $T$ by adding one new vertex $v$, and one new edge $uv$,
\item $X^+=X$, and 
\item $Y^+=Y$.
\end{itemize}
(Note that Operation 1 corresponds to Lemma \ref{lemma3}(1).)
\item {\bf Operation 2}

$(T^+,X^+,Y^+)$ arises by applying Operation 2 to $(T,X,Y)$ if 
\begin{itemize}
\item there is some vertex $u$ of $T$ with $u\not\in X$,
\item $T^+$ arises from $T$ by adding three new vertices $v$, $w_1$, and $w_2$, 
and three new edges $uv$, $vw_1$, and $vw_2$,
\item $X^+$ is either $X\cup \{ u,w_1,w_2\}$ or $X\cup \{ u,v,w_1,w_2\}$, and
\item $Y^+=Y\cup \{ u,v,w_1,w_2\}$.
\end{itemize}
(Note that Operation 2 corresponds to the special case of Lemma \ref{lemma2}(a)
where $k=\ell=2$, $W_1=\{ w_1\}$, and $W_2=\{ w_2\}$.)
\item {\bf Operation 3}

$(T^+,X^+,Y^+)$ arises by applying Operation 3 to $(T,X,Y)$ if 
\begin{itemize}
\item there is some vertex $u$ of $T$ with $u\not\in X$,
\item $T^+$ arises from $T$ by adding four new vertices $v$, $w_1$, $w_2$, and $w_3$, 
and four new edges $uv$, $vw_1$, $vw_2$, and $vw_3$,
\item $X^+$ is one of the four sets 
$X\cup \{ w_1,w_2,w_3\}$, 
$X\cup \{ u,w_1,w_2,w_3\}$, 
$X\cup \{ v,w_1,w_2,w_3\}$, or
$X\cup \{ u,v,w_1,w_2,w_3\}$, and
\item $Y^+=Y\cup \{ u,v,w_1,w_2,w_3\}$.
\end{itemize}
(Note that Operation 3 corresponds to the special case of Lemma \ref{lemma2}(b)
where $k=\ell=3$, $W_1=\{ w_1\}$, $W_2=\{ w_2\}$, and $W_3=\{ w_3\}$.)
\end{itemize}
The conditions that need to be satisfied 
in order to apply one of the last two extension operations 
are notationally more complicated.
For the description of these last two operations, 
we assume that $T$, $u$, $v$, $w_1$, $\ldots$, $w_k$, $W_1$, $\ldots$, $W_k$, and $\ell$ 
are as in the statement of Lemma \ref{lemma2}
and that either (a)(i) or (b)(i) is satisfied.
For an illustration see Figure \ref{fig3}.
Note that in the subforest $T[Y]$ of $T$ induced by $Y$,
the vertices $w_1,\ldots,w_k$ are endvertices,
that is, the vertex $v$ has at most one neighbor
that is not an endvertex.
\begin{itemize}
\item {\bf Operation 4}

$(T^+,X^+,Y^+)$ arises by applying Operation 4 to $(T,X,Y)$ if 
\begin{itemize}
\item $T^+$ arises from $T$ by adding one new vertex $v^*$, and one new edge $vv^*$, 
\item $X^+$ is either $X$ or $X\cup \{ v^*\}$, and
\item $Y^+=Y\cup \{ v^*\}$.
\end{itemize}
(Note that Operation 4 corresponds to Lemma \ref{lemma3}(2.1).)
\item {\bf Operation 5}

$(T^+,X^+,Y^+)$ arises by applying Operation 5 to $(T,X,Y)$ if 
\begin{itemize}
\item $T^+$ arises from $T$ by adding one new vertex $v^*$, and one new edge $w_iv^*$ for some $i\in [k]$, 
\item $X^+=X$, and
\item $Y^+=Y$.
\end{itemize}
(Note that Operation 5 corresponds to Lemma \ref{lemma3}(2.2).)
\end{itemize}
Let ${\cal T}$ be the set of triples $(T^+,X^+,Y^+)$ such that 
\begin{itemize}
\item either $(T^+,X^+,Y^+)\in \{ (K_1,\emptyset,\emptyset),(K_1,V(K_1),V(K_1))\}$,
where $K_1$ is the tree of order $1$,
\item or there is some triple $(T,X,Y)$ in ${\cal T}$
and some $i\in [5]$ such that 
$(T^+,X^+,Y^+)$ arises by applying Operation $i$ to $(T,X,Y)$.
\end{itemize}
The following is our main results, and yields a constructive characterization of ${\cal S}$.
\begin{theorem}\label{theorem1}
${\cal S}={\cal T}$.
\end{theorem}
{\it Proof:} Let $(T^+,X^+,Y^+)\in {\cal T}$. 
By induction on the order of $T^+$, 
we prove $(T^+,X^+,Y^+)\in {\cal S}$.
If $(T^+,X^+,Y^+)\in \{ (K_1,\emptyset,\emptyset),(K_1,V(K_1),V(K_1))\}$,
then $(T^+,X^+,Y^+)\in {\cal S}$ follows easily from the definition of ${\cal S}$.
Now, let $T^+$ have order at least $2$.
By the definition of ${\cal T}$, there is some $(T,X,Y)\in {\cal T}$ and some $i\in [5]$
such that $(T^+,X^+,Y^+)$ arises by applying Operation $i$ to $(T,X,Y)$.
By induction, $(T,X,Y)\in {\cal S}$.
If $i=1$, then Lemma \ref{lemma3}(1) implies $(T^+,X^+,Y^+)\in {\cal S}$.
If $i\in \{ 2,3\}$, then Lemma \ref{lemma2} implies $(T^+,X^+,Y^+)\in {\cal S}$.
Finally, 
if $i\in \{ 4,5\}$, then Lemma \ref{lemma3}(2) implies $(T^+,X^+,Y^+)\in {\cal S}$.
Altogether, we obtain ${\cal T}\subseteq {\cal S}$.

Now, let $(T,X,Y)\in {\cal S}$.
By induction on the order of $X$, 
we prove $(T,X,Y)\in {\cal T}$.
If $X=\emptyset$, 
then, by the definition of ${\cal S}$, 
we obtain that $Y=\emptyset$, 
and $(T,X,Y)$ arises from $(K_1,\emptyset,\emptyset)$ 
by some applications of Operation 1, 
which implies $(T,X,Y)\in {\cal T}$.
If $|X|=1$, then Lemma \ref{lemma1}(ii) implies that $(T,X,Y)=(K_1,V(K_1),V(K_1))$,
which implies $(T,X,Y)\in {\cal T}$.
Now let $|X|\geq 2$.
By Lemma \ref{lemma1}(ii), $|X|\geq 3$.
Let $P:z_0\ldots z_q$ be a longest path in $T$ such that $z_0,z_q\in X$.
Since $|X|\geq 3$ and $T$ is a tree, we have $q\geq 2$.
Let $u=z_2$, $v=z_1$, $w_1=z_0$, and $N_T(v)=\{ u,w_1,\ldots,w_k\}$.
Note that, by the choice of $P$, 
the vertices $u$, $v$, and $w_1,\ldots,w_k$ are as required in the statement of Lemma \ref{lemma2},
that is, $w_1$ belongs to $X$,
and, if $W_i$ is the vertex set of the component of $T-v$ that contains $w_i$,
then $W_i\cap X\subseteq \{ w_i\}$.
By Lemma \ref{lemma2}, there is some $(T',X',Y')\in {\cal S}$ such that 
either (a)(i) and (a)(ii),
or (b)(i) and (b)(ii)
are satisfied.
Since $|X'|<|X|$, we obtain, by induction, that $(T',X',Y')\in {\cal T}$.
If (a)(i) and (a)(ii) are satisfied, then $(T,X,Y)$ arises from $(T',X',Y')$ by 
\begin{itemize}
\item one application of Operation 2,
\item followed by some applications of Operation 4,
\item followed by some applications of Operation 5, 
\item followed by some applications of Operation 1.
\end{itemize}
If (b)(i) and (b)(ii) are satisfied, then $(T,X,Y)$ arises from $(T',X',Y')$ by 
\begin{itemize}
\item one application of Operation 3,
\item followed by some applications of Operation 4,
\item followed by some applications of Operation 5, 
\item followed by some applications of Operation 1.
\end{itemize}
By the definition of ${\cal T}$, this implies $(T,X,Y)\in {\cal T}$.
Altogether, we obtain ${\cal S}\subseteq {\cal T}$,
which completes the proof. $\Box$

\section{Conclusion}\label{section3}

The approach from Section \ref{section2} can be adapted 
to obtain some constructive characterization of the extremal trees for (\ref{e0}).
This naturally leads to a further refinement of the notion of a weak Roman dominating function.

Let $G$ be a graph, and let $X_0$ and $X_1$ be two disjoint subsets of the vertex set of $G$.
A {\it weak Roman dominating function for $(G,X_0,X_1)$}, a $(G,X_0,X_1)${\it-WRDF} for short, 
is a function $g:V(G)\to \{ 0,1,2\}$ 
such that every vertex $u$ in $X_0\cup X_1$ with $g(u)=0$ 
has a neighbor $v$ with $g(v)\geq 1$
such that $\{ x\in V(G):g_{v\to u}(x)\geq 1\}$ is $X_0$-dominating.
The {\it weak Roman domination number $\gamma_r(G,X_0,X_1)$} of $(G,X_0,X_1)$
is the minimum weight of a $(G,X_0,X_1)$-WRDF,
and a $(G,X_0,X_1)$-WRDF of weight $\gamma_r(G,X_0,X_1)$ is {\it minimum}.
Note that $f$ is a $(G,X)$-WRDF for some set $X$ of vertices of $G$ if and only if
$f$ is a $(G,X,\emptyset)$-WRDF.

Let ${\cal R}$ be the set of all $3$-tuples $((T_r,X_{r,0},X_{r,1}),(T_R,X_R),\delta_{R-r})$
with the following properties.
\begin{itemize}
\item $T_r$ is a tree, and $X_{r,0}$ and $X_{r,1}$ are disjoint sets of vertices of $T_r$.
\item $T_R$ is a tree, and $X_R$ is a set of vertices of $T_R$.
\item $\delta_{R-r}=\gamma_R(T_R,X_R)-\gamma_r(T_r,X_{r,0},X_{r,1})$.
\end{itemize}
A tree $T$ satisfies $\gamma_r(T)=\gamma_R(T)$
if and only if 
$((T,V(T),\emptyset),(T,V(T)),0)\in {\cal R}$.

There is a variant of Lemma \ref{lemma2} showing that 
$((T_r,X_{r,0},X_{r,1}),(T_R,X_R),\delta_{R-r})\in {\cal R}$
if and only if 
$((T'_r,X'_{r,0},X'_{r,1}),(T'_R,X'_R),\delta'_{R-r})\in {\cal R}$,
where 
$T_r'$ is a proper subtree of $T_r$
and 
$T_R'$ is a proper subtree of $T_R$.
Extracting the reductions encoded in this lemma similarly as in Section \ref{section2},
yields a constructive characterization of ${\cal R}$.
The drawback of this approach is that 
even if $T_r$ equals $T_R$, the tree $T_r'$ may be different from $T_R'$,
that is, in order to decide whether $\gamma_r(T)=\gamma_R(T)$ for some given tree $T$,
one has to generate/maintain two sequences of distinct subtrees  
$$T=T_r^{(0)}\supseteq T_r^{(1)}\supseteq T_r^{(2)}\ldots
\,\,\,\,\,\,\mbox{ and }\,\,\,\,\,\,
T=T_R^{(0)}\supseteq T_R^{(1)}\supseteq T_R^{(2)}\ldots$$
until a decision of possible,
because there is not necessarily always one reduction 
that works simultaneously for Roman domination as well as for weak Roman domination.
Therefore, such a constructive characterization essentially results in executing two separate reduction-based algorithms
that determine $\gamma_r(T)$ and $\gamma_R(T)$, and the decision amounts to a comparison of their results.
Not seeing much benefit in such an approach, we did not elaborate its details,
and leave it as an open problem to find a better constructive characterization of the extremal trees for (\ref{e0}).
Another interesting open problem is the complexity of deciding 
strong equality of the Roman domination number and the weak Roman domination number for general graphs.

\bigskip

\noindent {\bf Acknowledgment}  
J.D. Alvarado and S. Dantas were partially supported by FAPERJ, CNPq, and CAPES.


\begin{thebibliography}{}
\bibitem{cfmp} E.J. Cockayne, O. Favaron, C.M. Mynhardt, and J. Puech, A characterization of $(\gamma,i)$-trees, J. Graph Theory 34 (2000) 277-292.
\bibitem{chh} M. Chellali, T.W. Haynes, and S.T. Hedetniemi, Bounds on weak roman and 2-rainbow domination numbers, Discrete Appl. Math. 178 (2014) 27-32.
\bibitem{cr} M. Chellali and N.J. Rad, Strong equality between the Roman domination and independent Roman domination numbers in trees, Discuss. Math. Graph Theory 33 (2013) 337-346.
\bibitem{dghm} M. Dorfling, W. Goddard, M.A. Henning, and C.M. Mynhardt, Construction of trees and graphs with equal domination parameters, Discrete Math. 306 (2006) 2647-2654.
\bibitem{hh} J.H. Hattingh and M.A. Henning, Characterizations of trees with equal domination parameters, J. Graph Theory 34 (2000) 142-153.
\bibitem{hhs} T.W. Haynes, M.A. Henning, and P.J. Slater, Strong equality of domination parameters in trees, Discrete Math. 260 (2003) 77-87.
\bibitem{hs} T.W. Haynes and P.J. Slater, Paired-domination in graphs, Networks 32 (1998) 199-206.
\bibitem{h} M.A. Henning, A characterization of Roman trees, Discuss. Math. Graph Theory 22 (2002) 325-334.
\bibitem{hehe} M.A. Henning and S.T. Hedetniemi, Defending the Roman Empire - new strategy, Discrete Math. 266 (2003) 239-251.
\bibitem{rl} N.J. Rad and C.-H. Liu, Trees with strong equality between the Roman domination number and the unique response Roman domination number, Australas. J. Comb. 54 (2012) 133-140.
\bibitem{s} I. Stewart, Defend the Roman empire!, Sci. Am. 281 (1999) 136-139.
\end{thebibliography}
\end{document}